\documentclass[11pt, twoside]{article}
\usepackage{fancyhdr}
\usepackage{cite}
\usepackage[dvips, usenames]{color}
\usepackage{titletoc}
\usepackage{latexsym}
\usepackage{amsmath}
\usepackage{amssymb}
\usepackage{multicol}
\usepackage{graphics}
\usepackage{graphicx}
\usepackage{subfigure}
\usepackage{indentfirst}
\usepackage{epsfig}
\usepackage{mathrsfs}
\usepackage[
    driverfallback=pdftex, 
    colorlinks=true,
    citecolor=blue,
    linkcolor=blue,
    urlcolor=black,
    backref=page
]{hyperref}
\usepackage[nameinlink]{cleveref}
\usepackage[top=2.54cm, bottom=2.54cm, left=2.80cm, right=2.80cm]{geometry}
\def\q{\hspace*{\fill}$\Box$\medskip}

\newtheorem{theorem}{Theorem}
\newtheorem{proposition}[theorem]{Proposition}
\newtheorem{lemma}[theorem]{Lemma}

\newtheorem{conjecture}[theorem]{Conjecture}
\newtheorem{claim}{Claim}
\newtheorem{case}{Case}
\newtheorem{case2}{Case}[case]

\newtheorem{problem}[theorem]{Problem}

\newtheorem{definition}[theorem]{Definition}

\def\p{\noindent{\bf Proof. }}
\def\q{\hspace*{\fill}$\Box$\medskip}
\def\c{\hspace*{\fill}$\lozenge$\medskip}

\begin{document}
\title{\bf Spectral Radius Conditions for 3-Uniform Intersecting Families}
\author{Lusheng Fang$^{a}$, \, Guorong Gao$^{a,b}$\thanks {Email addresses: lusheng\_fang@126.com, grgao@fzu.edu.cn,
  anchang@fzu.edu.cn}, An Chang$^{a}$\\
\small$^a$Center for Discrete Mathematics and Theoretical Computer Science,\\
\small Fuzhou University, Fuzhou, Fujian, China\\
\small$^b$School of Mathematics and Statistics,
\small Fuzhou University, Fuzhou, Fujian, China\\
}
\date{}
\maketitle

\begin{abstract}

{ Let $M_k$ denote a matching of size $k$. The classical Erd\H{o}s matching conjecture asks for the maximum number of edges of an intersecting $r$-graph without $M_k$.  The csae for $k=2$, which is known as intersecting $r$-graph, is established by Erd\H{o}s, Ko and Rado. Hilton and Milner further determine the maximum number of edges of a non-trivial intersecting $r$-graph, where the intersecting $r$-graph $H$ is called non-trivial if $\cap_{e\in E(H)}e=\emptyset$. 
In this paper, we investigate the spectral analogues of the hpergraph matching problems and intersecting family problems. More precisely, for sufficiently large $n$, we determine respectively the maximum spectral
radius of $M_{k+1}$-free and non-trivial intersecting $3$-graphs on $n$ vertices, and characterize the extremal hypergraphs. 
}
\vspace{2mm}
\vspace{2mm}

\noindent{\bf Keywords}\,:~Hypergraph, \space\space Spectral radius,\space\space Matching
\vspace{2mm}
\vspace{2mm}

\noindent{\bf AMS subject classiflcations}\,:~05C35, \space 05C65
\end{abstract}

\vskip 0.2in
\section{Introduction}

An $n$-vertex $r$-uniform hypergraph (or $r$-graph for short) $H = (V(H),E(H))$ consists of a vertex set $V(H)=\{1,2,...,n\}$ and an edge set $E(H)$, where $E(H)$ is a set of $r$-element subsets of $V(H)$. As usual, a graph is a $2$-uniform hypergraph. Throughout this paper, we default to $r\geq 2$.

 For an $r$-graph $H$, the spectral radius of $H$, denoted by $\rho(H)$, is defined as the maximum absolute value of the eigenvalues of the adjacency tensor $\mathcal{A}(H)$ (see the definitions in Section \ref{spe}). 
 In \cite{Bru86}, Brualdi and Solheid  introduced a notable class of extremal problems in spectral graph theory, which can be stated as follows: What is the maximum spectral radius in a graph $G$ if the graph $G$ belongs to  a specified class of graphs. 
 Since then, numerous graph classes have been studied, such as  given chromatic index\cite{Feng16}, given diameter \cite{Hansen08} and so on. 
 
 In this paper, we focus our attention on the following problem which is a natural extension of graph case.
 \begin{problem}\label{BSP}
   What is the maximum spectral radius in an $r$-graph $H$ if the $r$-graph $H$ belongs to  a specified class of $r$-graphs?
 \end{problem}
 
For an $r$-graph $H$ and a family of $r$-graphs $\mathcal{F}$, we say that $H$ is $\mathcal{F}$-free if $H$ does not contain any member of $\mathcal{F}$ as a subgraph. 
Determining the maximum spectral radius among all $\mathcal{F}$-free $r$-graphs has become one of the central topics in Problem \ref{BSP}, usually referred to as the spectral Tur\'an problem.  

Most existing studies on spectral Tur\'an problems concern nondegenerate hypergraphs, that is, hypergraphs $\mathcal{F}$ with positive Tur\'an density $$\pi(\mathcal{F})=\lim_{n\to\infty}\frac{\mathrm{ex}_r(n,\mathcal{F})}
{\binom{n}{r}}>0,$$
where $ex_r(n, \mathcal{F})$ is the maximum possible number of edges of an $n$-vertex $\mathcal{F}$-free $r$-graph.
 Early progress was made by Keevash, Lenz and Mubayi \cite{Kee14}, who determined the maximum $p$-spectral radius of Fano-free $3$-graphs when $n$ is sufficiently large.  Extending the spectral Mantel's theorem to hypergraphs, Ni, Liu and Kang \cite{Ni22}  determined the maximum $p$-spectral radius of cancellative 3-graphs, and characterized the extremal hypergrap.  Given a graph $F$, the $r$-expansion of $F$, denoted by $F^+$, is the $r$-graph obtained from $F$ by enlarging each edge of $F$ with $r-2$ new vertices disjoint from $V(F)$ such that distinct edges of $F$ are enlarged by distinct vertices. Gao, Chang and Hou \cite{Gao22} determined the maximum spectral radius of $K_{r+1}^{+}$-free linear $r$-graphs, and She, Fan, Kang and Hou \cite{She22} extended their result to $F^{+}$-free linear $r$-graphs, where $F$ is a color-critical graph.  For a nondegenerate $r$-graph $F$, 
spectral stability results,
namely that every near-extremal (with respect to spectral radius) $n$-vertex $F$-free graph is structurally close to an extremal graph,
are powerful tools in the study of the spectral Tur\'an problem for $F$. For more details on the spectral stability results, one can see \cite{Zheng24,Liu24, Ni26, Zheng25}.

In contrast, considerably less is known about spectral Tur\'an problems for degenerate hypergraphs, namely hypergraphs $\mathcal{F}$ with $\pi(\mathcal{F})=0$.
Recently, Fang, Gao, Chang and Hou \cite{FangE25} established several sharp bounds on the spectral radius of $K_{s,t}^{+}$-free linear $r$-graphs.
Given a graph $F$, a hypergraph $H$ is called a Berge-$F$ if there exists a bijection $f:E(F)\rightarrow E(H)$ such that $e\subseteq f(e)$ for every $e\in E(F)$. The same group of authors \cite{FangD25} subsequently  determined the spectral extremal problem for Berge-$P_k$-free $3$-graphs for  sufficiently large $n$ and $k\geq 8$ and characterized the extremal hypergraphs. For more results on spectral Tur\'an problems for degenerate hypergraphs, we refer the reader to
\cite{Wang25,Zhu24,Zhou24,Zhou25, Hou21,She24}.

Among various degenerate hypergraphs, 
determining the maximum number of edges in a hypergraph with bounded matching number has been one of the most influential topics in extremal combinatorics.
A matching of size $k$ in an $r$-graph, denoted by $M_k$,  is a collection of $k$ pairwise disjoint edges. In particular,
an $r$-graph  is intersecting if it is $M_2$-free.
A central extremal problem for $r$-graphs is to determine the maximum number of edges in an $M_{k+1}$-free $r$-graph on $n$ vertices. The case $k=1$ corresponds to the
classical Erd\H{o}s-Ko-Rado theorem \cite{Erd61}, which is stated as follows.

 \begin{theorem} [Erd\H{o}s-Ko-Rado Theorem \cite{Erd61}]
 If $H$ is an  intersecting $r$-graph on $n\geq 2r$ vertices, then 
 $$|E(H)|\leq \binom{n-1}{r-1}.$$
\end{theorem}
Since intersecting hypergraphs are $M_2$-free, a natural generalization of the Erd\H{o}s-Ko-Rado theorem is to consider $r$-graphs that are $M_{k+1}$-free, where $k>1$. Motivated by this direction, Erd\H{o}s proposed the following classical conjecture and related result.

\begin{conjecture} [Erd\H{o}s Matching Conjecture \cite{Erd65}]
Let $n\geq (k+1)r$.  If $H$ be an $M_{k+1}$-free $r$-graph on $n$ vertices, then 
   $$|E(H)|\leq \max\{|\mathcal{A}^r_1|,|\mathcal{A}^r_r|\}=\max\left\{\binom{n}{r}-\binom{n-k}{r}, \binom{rk+r-1}{r}\right\},$$
where $\mathcal{A}^r_i=\left\{e\in \binom{[n]}{r}:|e\cap[ik+i-1]\geq i\right\}$ for $1\leq i\leq r$.
\end{conjecture}
\begin{theorem} [\cite{Erd65}]\label{Erd}
Let $n\geq n_0(r,k)$. If $H$ is an $M_{k+1}$-free $r$-graph on $n$ vertices, then 
   $$|E(H)|\leq |\mathcal{A}^r_1|.$$
\end{theorem}
The Erd\H{o}s Matching Conjecture is trivial for $r=1$, and has been completely resolved for $r=2$ and $r=3$ by Erd\H{o}s, Gallai \cite{Erd59} and Frankl \cite{Fra17}, respectively. For more details, one can see \cite{Frankl26, Hou26}.

Besides the Erd\H{o}s Matching Conjecture, another classical problem concerning intersecting hypergraphs is to determine the maximum size of a non-trivial intersecting hypergraph. Note that an intersecting $r$-graph $H$ is called non-trivial if $\cap_{e\in E(H)}e=\emptyset$, and trivial otherwise. Let $$\mathcal{F}^r_1=\mathcal{B}\cup\{[2,r+1]\}$$
where $\mathcal{B}=\left\{e\in \binom{[n]}{r}:1\in e, e\cap[2,r+1]\neq \emptyset \right\}$ and 
$$\mathcal{H}_1=\left\{e\in \binom{[n]}{3}:|e\cap\{1,2,3\}|\geq 2\right\}.$$
Clearly, $\mathcal{F}^r_1$ and $\mathcal{H}_1$ are non-trivial hypergraphs.
In 1967, Hilton and Milner \cite{Hil67} determined the
maximum number of edges of a non-trivial intersecting $r$-graph on $n$ vertices.
 \begin{theorem} [Hilton-Milner Theorem \cite{Hil67}]\label{HMT}
 Let $r\geq 2$ and $n\geq 2r$. If $H$ is a non-trivial intersecting $r$-graph on $n$ vertices, then 
 \begin{equation*}
\begin{aligned}
|E(H)|\leq \binom{n-1}{r-1}-\binom{n-r-1}{r-1}+1.
\end{aligned}
\end{equation*}
Moreover, equality holds if and only if $H$ is $\mathcal{F}^r_1$, or, in the cases $r=3$, $H$ is $\mathcal{H}_1$.
\end{theorem}

It is therefore natural to ask whether the extremal hypergraphs in Theorem \ref{Erd} and the Hilton--Milner theorem also maximize the spectral radius among all hypergraphs satisfying the corresponding constraints. Motivated by this question, we investigate spectral analogues of these classical results. Our first result establishes a spectral version of Theorem \ref{Erd} for $3$-graphs, while our second result provides a spectral analogue of the Hilton--Milner theorem for $3$-graphs.
\begin{theorem}\label{ErdT}
Let $k\geq 1$ be an integer and $k=o(n)$.
For sufficiently large $n$, if $H$ is a $M_{k+1}$-free $3$-graph on $n$ vertices, then
$$\rho(H)\leq \rho(\mathcal{A}^3_1),$$
with equality if $H=\mathcal{A}^3_1$.
\end{theorem}
\begin{theorem}\label{ErdTP}
Let $k\geq 1$ be an integer and $k=\theta(\frac{n}{3})$.
For sufficiently large $n$, if $H$ is a $M_{k+1}$-free $3$-graph on $n$ vertices, then
$$\rho(H)\leq \rho(\mathcal{A}^3_3),$$
with equality if $H=\mathcal{A}^3_3$.
\end{theorem}
\begin{theorem}\label{HMTT}
For sufficiently large $n$,
if $H$ is a non-trivial intersecting $3$-graph on $n$ vertices, then
$$\rho(H)\leq \rho(\mathcal{H}_1),$$
with equality if $H=\mathcal{H}_1$.
\end{theorem}

The rest of this paper is organized as follows. In the next section, we make a brief introduction on the adjacency tensor of a uniform hypergraph and its eigenvalues, spectral radius, etc.
Meanwhile, we also introduce the shifting operation of hypergraphs
which plays an important role in the subsequent proof process. Finally, we prove our main results in Section \ref{Proof}.

\section{Preliminaries}

\subsection{Eigenvalues of tensors}\label{spe}

In 2005, Qi \cite{Qi05} and Lim \cite{Lim05} independently introduced the concept of tensor eigenvalues and the spectra of tensors. An order $r$ dimension $n$ real tensor $\mathcal{T}=(\mathcal{T}_{i_1\cdots i_r})$ consists of $n^{r}$ real entries $\mathcal{T}_{i_1\cdots i_r}$  for all $i_1, i_2, \cdots, i_r\in  [n]$, where $[n]=\{1,2,\ldots,n\}$. Evidently, a vector of dimension $n$ is a tensor of order 1 and matrix is a tensor of order 2. A tensor $\mathcal{T}=(\mathcal{T}_{i_1\cdots i_r})$ is called symmetric if its entries $\mathcal{T}_{i_1\cdots i_r}$'s are invariant under any permutation of the indices $i_1, i_2, \cdots, i_r$. Given a vector $\mathbf{x}=(x_1,x_2,...,x_n)^{\intercal} \in \mathbb{R}^n$, we adopt the following notation: $\mathcal{T}\mathbf{x}^r$ is a real number and $\mathcal{T}\mathbf{x}^{r-1}$ is an $n$-dimensional vector, where $\mathcal{T}\mathbf{x}^r$ and the $i$th component of $\mathcal{T}\mathbf{x}^{r-1}$ are given by:
\begin{equation*}
\begin{aligned}
\mathcal{T}\mathbf{x}^r = \sum_{i_1, i_2, \cdots, i_r \in [n]}\mathcal{T}_{i_1i_2\cdots i_r}x_{i_1}x_{i_2}\cdots x_{i_r}.
\end{aligned}
\end{equation*}
\begin{equation*}
\begin{aligned}
(\mathcal{T}\mathbf{x}^{r-1})_i = \sum_{ i_2, \cdots, i_r \in [n]}\mathcal{T}_{ii_2\cdots i_r}x_{i_2}\cdots x_{i_r}.
\end{aligned}
\end{equation*}
If there \textcolor{red}{exist} $\lambda \in \mathbb{C}$ and a nonzero vector $\mathbf{x} \in \mathbb{C}^n$ satisfying
\begin{equation*}
\begin{aligned}
\mathcal{T}\mathbf{x}^{r-1}=\lambda \mathbf{x}^{[r-1]},
\end{aligned}
\end{equation*}
then $\lambda$ is called an eigenvalue of $\mathcal{T}$ and $\mathbf{x}$ is its corresponding eigenvector, where $\mathbf{x}^{[r-1]}=(x_{1}^{r-1}, x_{2}^{r-1}, \cdots,  x_{n}^{r-1})^T\in \mathbb{C}^n\setminus\{0\}$. If $\mathbf{x}$ is a real eigenvector of $\mathcal{T}$,
surely the corresponding eigenvalue $\lambda$ is real. In this case, $\lambda$ is called an $H$-eigenvalue and $\mathbf{x}$ is called an
$H$-eigenvector associated with $\lambda$. Furthermore, if $\mathbf{x}$ is nonnegative and real, we say $\lambda$ is an $H^{+}$-eigenvalue of $\mathcal{T}$. If $\mathbf{x}$ is positive and real, $\lambda$ is said to be an $H^{++}$-eigenvalue of $\mathcal{T}$. Throughout this paper, we will refer to H-eigenvalues and H-eigenvectors as eigenvalues and eigenvectors, or use both interchangeably. The maximal absolute value of the eigenvalues of $\mathcal{T}$ is called the \textbf{spectral radius} of $\mathcal{T}$, denoted by $\rho(\mathcal{T})$.

A tensor $\mathcal{T}=(\mathcal{T}_{i_1\cdots i_r})$ of order $r$ dimension $n$ is called reducible if there exists a
nonempty proper index subset $I\subset \{1, \cdots, n\}$ such that
$\mathcal{T}_{i_1\cdots i_r}= 0$  $\forall i_{1}\in I$,  $\forall i_{2}, . . . , i_{m} \notin I$.
If $\mathcal{T}$ is not reducible, then we call $\mathcal{T}$ irreducible. A tensor with all nonnegative entries is called a nonnegative tensor. It is well known that the Perron-Frobenius theorem for nonnegative matrices plays a crucial role in the study of spectral graph theory. As an extension of matrices, the Perron-Frobenius theorem for nonnegative tensors has also been established, see \cite{Cha08,Fri13,Yan10}.
\begin{theorem}[Perron-Frobenius theorem for nonnegative tensors]\label{Perron-Frobenius theorem for nonnegative tensors}\mbox{}\par
\noindent(1) (Yang and Yang, 2010 \cite{Yan10}). If $\mathcal{T}$ is a nonnegative tensor of order $r$ and dimension $n$,
then $\rho(\mathcal{T})$ is an $H^{+}$-eigenvalue of $\mathcal{T}$.\\
(2) (Friedland, Gaubert and Han, 2011 \cite{Fri13}). If furthermore $\mathcal{T}$ is weakly irreducible, then
$\rho(\mathcal{T})$ is the unique $H^{++}$-eigenvalue of $\mathcal{T}$, with the unique eigenvector $x\in R_{++}^{n}$, up
to a positive scaling coefficient.\\
(3) (Chang, Pearson and Zhang, 2008 \cite{Cha08}). If moreover $\mathcal{T}$ is irreducible, then $\rho(\mathcal{T})$ is the
unique $H^{+}$-eigenvalue of $\mathcal{T}$, with the unique eigenvector $x\in R_{+}^{n}$, up to a positive
scaling coefficient.
\end{theorem}

In 2012, Cooper and Dutle \cite{Coo12} defined the adjacency tensor $\mathcal{A}(H)$ of an $r$-graph $H$ with $n$ vertices as an $r$-th order $n$  dimensional symmetric tensor $\mathcal{A}(H)=\mathcal{A}_{i_1\cdots i_r}$, where
\begin{equation*}
\mathcal{A}_{i_1\cdots i_r}=\left\{
\begin{array}{ll}
\frac{1}{(r-1)!}&  \;\ if\ \{i_1\cdots i_r\} \in E,\\
0 &\;\ otherwise.\\
\end{array}\right.
\end{equation*}

For an $r$-graph $H$, the spectral radius of $H$ is defined as that of the adjacency tensor $\mathcal{A}(H)$, denoted by $\rho(H)$. Since the adjacency tensor $\mathcal{A}(H)$ of $H$ is a nonnegative tensor, its spectral
radius $\rho(H)$ is an $H^{+}$-eigenvalue of $\mathcal{A}(H)$. In addition, Qi \cite{Qi05} showed that $\rho(H)$ is also the optimal value of the following maximization
\begin{equation} \label{Optimal}
\begin{aligned}
\rho(H)=\max_{\|x\|_r=1}\mathbf{x}^T\mathcal{A}(H)\mathbf{x}^{r-1},
\end{aligned}
\end{equation}
and $\rho(H)=\mathcal{A}(H)\mathbf{x}^{r}$ for the optimal vector $\mathbf{x}$.  Specifically, Friedland et al. \cite{Fri13}  proved that a uniform hypergraph $H$ is connected if and only if its adjacency tensor $\mathcal{A}(H)$ is weakly irreducible.
By the above Perron-Frobenius theorem, if $H$ is connected, then the eigenvector corresponding to the spectral radius $\rho(H)$ can be chosen to be strictly positive, which is known as the principal eigenvector. 

We denote by $\mathbb{S}_{r,+}^{n-1}$ the set of all nonnegative real vectors $\mathbf{x}\in \mathbb{R}^{n}$ with $\|\mathbf{x}\|_r=1$. Let $\mathbf{x}$ be an eigenvector of $\mathcal{A(H)}$ and $U\subseteq V(H)$. We define $\mathbf{x}(U)=\prod_{v_i\in U}x_i$. Clearly, for a vector $\mathbf{x}$ and $v\in V(H)$, we have
\begin{equation}\label{r eigenvector equation}
\left(\mathcal{A(H)}\mathbf{x}^{r-1}\right)_v=\sum_{e\in  E(H), v\in e}\mathbf{x}(e\backslash\{v\}).
\end{equation}

\subsection{The shifting operation}

The \textbf{shifting operation} is a very useful combinatorial tool  for studying  extremal families with certain intersection properties, which was first introduced by Frankl in \cite{Fra87}. Subsequently, Keough and Radcliffe \cite{Keo22} introduced the following shifting operation $S_{i\rightarrow j}$ of a hypergraph $H$, where $i,j\in V(H)$.

\begin{definition} [\cite{Keo22}]\label{Shifting}
Given a set $A\subseteq [n]$ and $i,j\in [n]$ such that $A\cap \{i,j\}=i$, define $A_{i\rightarrow j}=(A\backslash\{i\})\cup \{j\}.$ Consider a hypergraph $H$ with vertex set $[n]$ and edge set $E$. For $1\leq j <i\leq n$, define
$H_{i\rightarrow j}$  to be the $r$-graph with vertex set $[n]$ and edge set $\{S_{i\rightarrow j}(e):e\in E\}\cup \{e:e,S_{i\rightarrow j}(e) \in E\}$, where for each $e\in E,$

\begin{equation*}
S_{i\rightarrow j}(e)=
\begin{cases}
e_{i\rightarrow j} & if\  e\cap \{i,j\}=\{i\},\\
e & otherwise.
\end{cases}
\end{equation*}
\end{definition}

\begin{lemma}\label{Matching}
  Let $H$ be a $M_{k+1}$-free $r$-graph. Then, for any $u,v\in V(H)$, $H_{u\rightarrow v}$ is also $M_{k+1}$-free.
\end{lemma}
\p
We prove this by contradiction. Suppose that there exist $k+1$  edges $e_1, e_2,\ldots,$ $e_{k+1}$ in $S_{u\rightarrow v}(H)$ such that they are pairwise disjoint. Clearly, it is impossible to $e_1, e_2, \ldots, e_{k+1}\in E(H).$ Thus, there exists an edge $e_t$ such that $v\in e_t, u\notin e_t$ and $e_t\notin E(H)$, where $t\in [k+1]$. Let $e'_t=(e_t\backslash\{v\})\cup\{u\}$. Then, 
by the definition \ref{Shifting} of shifting, it is easy to see that $e'_t\in E(H)$. If $u \notin \bigcup_{i=1}^{k+1}V(e_i)$, we have that $e_1, \cdots, e_{t-1}, e'_t,e_{t+1},\cdots, e_{k+1}$  are pairwise disjoint in $H$, which is a contradiction. 
Hence,  there exists an edge $e_l$ such that $u\in e_l$ and $v\notin e_t$, where $l\in [k+1]$.  Let $e'_l=(e_l\backslash\{u\})\cup\{v\}$. Then, 
by the definition \ref{Shifting} of shifting, $e'_l\in E(H)$. Without loss of generality, $l>t$.  Therefore, 
$$e_1, \cdots, e_{t-1}, e'_t, e_{t+1}\cdots,e_{l-1}, e'_l, e_{l+1},\cdots, e_{k+1}$$
are pairwise disjoint in $H$, which is a contradiction.
\q

%\begin{proposition} \label{pro1}
%$\alpha'(S_{i\rightarrow j}(H)) \leq \alpha'(H)$
%\end{proposition}
%\p
%Suppose $e_1, e_2, \cdots, e_l$ are independent edges in $S_{i\rightarrow j}(H)$. If $j \notin \bigcup_{i=1}^lV(e_i)$, then $e_1, e_2, \cdots, e_l$  are independent edges in $H$. Meanwhile, if $j\in \bigcup_{i=1}^lV(e_i)$ and $j \in e_t \in E(H)$ in $S_{i\rightarrow j}(H)$ where $t \in \{1, 2, \cdots, l\}$, then $e_1, e_2, \cdots, e_l$  are independent edges in $H$. If $j \in \bigcup_{i=1}^lV(e_i)$, $j \in e_t \notin E(H)$ and $i \notin \bigcup_{i=1}^lV(e_i)$, then $e_1, \cdots, e_{t-1}, (e_{t} \backslash \{j\})\bigcup\{i\},e_{t+1},\cdots,e_l$  are independent edges in $H$.   If $j \in \bigcup_{i=1}^lV(e_i)$, $j \in e_t \notin E(H)$ and $i \in e_k \in \{e_1, e_2, \cdots, e_l\}$, then $e_1, \cdots, e_{t-1}, (e_{t} \backslash \{j\})\bigcup\{i\},e_{t+1},\cdots,e_{k-1}, (e_{k} \backslash \{i\})\bigcup\{j\},e_{k+1},\cdots, e_l$  are independent edges in $H$, where $j\in e_j$ and $i\in e_i$ in $S_{i\rightarrow j}(H)$.
%\q

\begin{lemma}[\cite{FangD25}]\label{spectral}
Let $H=(V,E)$ be an $r$-graph and $\{i,j\}\subseteq V$. Then
\begin{equation*}
\begin{aligned}
\rho(H_{i\rightarrow j})\geq\rho(H).
\end{aligned}
\end{equation*}
\end{lemma}

A hypergraph $H=([n],E)$ is $\mathbf{shifted}$ if and only if $H_{i\rightarrow j}=H$ for all $1\leq j<i\leq n$.

Given $A=\{a_1, a_2, \cdots, a_r\}\in \binom{[n]}{r}$ and $B=\{b_1, b_2, \cdots, b_r\}\in \binom{[n]}{r}$ we let $A\preceq B$ if $a_i \leq b_i$. In shifted hypergraph, we can obtain the following proposition.
\begin{proposition}[\cite{Keo22}]\label{shiftedP}
Let $H$ be a shifted hypergraph. If $e_1 \in E(H)$ and $e_2 \preceq e_1$, then $e_2 \in E(H)$.
\end{proposition}

\section{Spectral extremal problems}\label{Proof}
In this section, we present the proofs of our main results. 

Given two vertices $u$ and $v$, we say that $u$ and $v$ are \textbf{equivalent} in $H$, in writing $u\thicksim v$, if
transposing $u$ and $v$ and leaving the remaining vertices intact, we get an automorphism of $H$.

\begin{lemma} [\cite{Nik14}] \label{lm_Liu}
Let $H$ be an $r$-graph on $n$ vertices and $u\thicksim v$. If $\mathbf{x}\in \mathbb{S}^{n-1} $ is an eigenvector to $\rho(H)$, then $x_u=x_v$.
\end{lemma}

\begin{lemma}[\cite{Kee14}]\label{Hou}
  Let $H$ be an $r$-graph with $m$ edges. Then
  $$\rho(H)\leq (rm)^{\frac{r-1}{r}}.$$
\end{lemma}

\subsection{Proof of the Theorem \ref{ErdT}}
\p
Let $H$ be a $M_{k+1}$-free $3$-graph on vertex set $[n]$ with the maximum spectral radius. Let $H'$ be a shifted hypergraph obtained by applying a number of shifting operations to $H$. By Lemma \ref{Matching}, $H'$ is also $M_{k+1}$-free. By the choice of $H$ and Lemma \ref{spectral}, it follows that $\rho(H)=\rho(H')$. 
Let $e_1=\{1, 3k+2, 3k+3\}, e_2=\{2, 3k, 3k+1\}, \ldots, e_i=\{i, 3k+4-2i, 3k+5-2i\}, \ldots, e_{k+1}=\{k+1, k+2, k+3\}$. Clearly, for any $i,j\in [k+1]$, $e_i\cap e_j=\emptyset$.
Thus, it is impossible that $e_1, e_2,\ldots, e_{k+1}\in E(H')$. 

For any $i\in [k+1]$, let $H_i$ be a shifted $3$-graph on vertex set $[n]$  with the maximum number of edges that does not contain $e_i$. Then, $$\rho(H')\leq \max\{\rho(H_1), \rho(H_2), \ldots, \rho(H_{k+1})\}.$$ 
By Proposition \ref{shiftedP}, for any $e=\{a,b,c\}\in E(H_i)$, $a\leq i-1$ or $b\leq 3k+3-2i$ or $c\leq 3k+4-2i$. Let $V_1^i=\{1,2,\ldots, i-1\}$, $V_2^i=\{i,i+1,\ldots, 3k+3-2i\}$ and $V_3^i=\{3k+4-2i, 3k+5-2i, \ldots, n\}$. 
Note that as a supplementary assumption, we set $V_1^i=\emptyset$ when $i=1$.
Let $E_1=\left\{e\in \binom{[n]}{3}: |e\cap \left(V_1^i\cup V_2^i\right)|\geq 2\right\}$ and $E_2=\left\{e\in \binom{[n]}{3}: |e\cap V_1^i|= 1, |e\cap V_3^i|= 2\right\}$.
Then, $E(H_i)=E_1\cup E_2.$
Obviously, $H_{k+1}\cong \mathcal{A}^3_1$ and $\mathcal{A}^3_1$ is $M_{k+1}$-free. So, $\rho(H')\geq \rho(H_{k+1})$. To finish the proof, we present the following claim.

\begin{claim}\label{spectral comparison}
  For large enough $n$ and $1\leq i\leq k$, we show that $\rho(H_{i})< \rho(H_{i+1})$.
\end{claim}
\p
For any $i,j\in V_l^i$, it is easy to see that $i$ and $j$ are equivalent in $H_i$, where $l=1,2,3$.
Let $\mathbf{x}\in \mathbb{S}_{3,+}^{n-1}$ be an eigenvector to $\rho(H_i)$. By Lemma \ref{lm_Liu}, we denote $x_l=x_v$ for $v\in V_l^i$, where $l=1,2,3$.
Let $$E_3=\left\{e\in \binom{[n]}{3}:i\in e, |e\cap V_3^i|=2\right\},$$ 
$$E_4=\left\{e\in \binom{[n]}{3}: i\notin e, \{3k+2-2i, 3k+3-2i\}\subseteq e, |e\cap V_2^i|=3\right\},$$ 
\begin{equation*}
\begin{split}
E_5 = \bigg\{ e \in \binom{[n]}{3}: &\ i \notin e,\ |e \cap \{3k+2-2i, 3k+3-2i\}| = 1, \\
&\ |e \cap V_2^i| = 2,\ |e \cap V_3^i| = 1 \bigg\}.
\end{split}
\end{equation*}
and
\begin{equation*}
\begin{split}
E_6 = \bigg\{ e \in \binom{[n]}{3}: &\{3k+2-2i, 3k+3-2i\}\subseteq e, |e \cap V_3^i| = 1 \bigg\}.
\end{split}
\end{equation*}
Then, $E(H_{i+1})=\left(E(H_i)\cup E_3\right)\backslash (E_4\cup E_5\cup E_6).$ Hence, it follows from that
\begin{equation}\label{diff2}
\begin{aligned}
\rho(H_{i+1})-\rho(H_i)&\geq 3\left(\sum_{e\in E_3}\mathbf{x}(e)-\sum_{e\in E_4\cup E_5\cup E_6}\mathbf{x}(e)\right)\\
&=3\binom{n-3k-3+2i}{2}x_2x_3^{2}\\
&-3(3k-3i+1)x_2^3\\
&-3\big(2(n-3k+2i-3)+1\big)(3k-3i+1)x_2^2x_3.
\end{aligned}
\end{equation}
Let $\rho(H_i)=cn^{\alpha}+o(n^{\alpha})$ and $k=tn^{\beta}+o(n^{\beta})$, where $c$ and $t$ are constants. For $l=1,2,3$, let $x_l=c_ln^{\alpha_l}+o(n^{\alpha_l})$, where $c_l$ is a constant. 
By the eigenvector equation, we have
\begin{equation}\label{eigenvector equation5}
\begin{aligned}
\rho(H_{i})x_3^{2}&=\binom{i-1}{2}x_1^{2}+(3k-3i+4)(i-1)x_1x_2+\binom{3k-3i+4}{2}x_2^{2}\\
&+(i-1)(n-3k-4+2i)x_1x_3.
\end{aligned}
\end{equation}
By (\ref{eigenvector equation5}), it is easy to see that
\begin{equation}\label{case1}
c'_1n^{\alpha+2\alpha_3}+o(n^{\alpha+2\alpha_3})\geq c'_2n^{2\alpha_2+2\beta}+o(n^{2\alpha_2+2\beta}),
\end{equation}
where $c'_1$ and $c'_2$ are  constant. Note that $|E(H_{i})|=O(n^{2+\beta})$. By Lemma \ref{Hou}, we have $\rho(H_{i})=O(n^{\frac{2(2+\beta)}{3}})$. Hence, $\alpha\leq \frac{2(2+\beta)}{3}$. For large enough $n$, by (\ref{case1}), we have that 
$\frac{2(2+\beta)}{3}+2\alpha_3\geq 2(\alpha_2+\beta).$
Then, 
\begin{equation}\label{Mulkey2}
\frac{2+\beta}{3}+\alpha_3\geq \alpha_2+\beta.
\end{equation}
Since $\beta<1$, combining (\ref{diff2}) and (\ref{Mulkey2}),  we obtain
\begin{equation*}
\begin{aligned}
\rho(H_{i+1})-\rho(H_i)&>\theta(n^{2+\alpha_2+2\alpha_3})-\theta(n^{1+2\alpha_2+\alpha_3+\beta})-\theta(n^{3\alpha_2})\\
&>0.
\end{aligned}
\end{equation*}
\c

Recall that $\rho(H')\leq \max\{\rho(H_1), \rho(H_2), \ldots, \rho(H_{k+1})\}$ and $\rho(H')\geq \rho(H_{k+1})$. Thus, by Claim \ref{spectral comparison}, $$\rho(H')=\rho(H_{k+1})=\rho(\mathcal{A}^3_1).$$ This completes the proof.
\q
\subsection{Proof of the Theorem \ref{ErdTP}}
\p
Let $H$ be a $M_{k+1}$-free $3$-graph on vertex set $[n]$ with the maximum spectral radius. Let $H'$ be a shifted hypergraph obtained by applying a number of shifting operations to $H$. By Lemma \ref{Matching}, $H'$ is also $M_{k+1}$-free. By the choice of $H$ and Lemma \ref{spectral}, it follows that $\rho(H)=\rho(H')$. 

Let $e_1=\{1, 2, 3k+3\}, e_2=\{3, 4, 3k+2\}, \ldots, e_i=\{2i-1, 2i, 3k+3-i+1\}, \ldots, e_{k+1}=\{2k+1, 2(k+1), 2k+3\}$. Clearly, for any $i,j\in [k+1]$, $e_i\cap e_j=\emptyset$.
Thus, it is impossible that $e_1, e_2,\ldots, e_{k+1}\in E(H)$. 
For any $i\in [k+1]$, let $H_i$ be a shifted $3$-graph on vertex set $[n]$  with the maximum number of edges that does not contain $e_i$. Then, $$\rho(H)\leq \max\{\rho(H_1), \rho(H_2), \ldots, \rho(H_{k+1})\}.$$ 
 By Proposition \ref{shiftedP}, for any $e=\{a,b,c\}\in E(H_i)$, $a\leq 2i-2$ or $b\leq 2i-1$ or $c\leq 3k+3-i$. Let $V_1^i=\{1,2,\ldots, 2i-2\}$, $V_2^i=\{2i-1,2i,\ldots, 3k+3-i\}$ and $V_3^i=\{3k+4-i, 3k+5-i, \ldots, n\}$. 
Note that as a supplementary assumption, we set $V_1^i=\emptyset$ when $i=1$.
Then, $$E(H_i)=\left\{e\in \binom{[n]}{3}: |e\cap V_1^i|\geq 1\right\}\cup \left\{e\in \binom{V_2^i}{3}\right\}.$$
Obviously, $H_1=K^3_{3k+2}\cup (n-3k-2)K_1=\mathcal{A}^3_3$.  Clearly, $H_1$ is $M_{k+1}$-free. Hence, $\rho(H')\geq \rho(H_{1})$. To finish the proof, we present the following claim.
%To complete the proof, for sufficiently large $n$ and $2\leq i\leq k+1$, we show that $\rho(H_{i})< \rho(H_{1})$. To complete the proof, for sufficiently large $n$ and $2\leq i\leq k+1$, we show that $\rho(H_{i})< \rho(H_{1})$.

\begin{claim}\label{spectral comparison2}
  For large enough $n$ and $1\leq i\leq k$, we show that $\rho(H_{i})< \rho(H_{i+1})$.
\end{claim}
\p
For any $i,j\in V_l^i$, it is clear that $i$ and $j$ are equivalent in $H_i$, where $l=1,2,3$.
Let $\mathbf{x}\in \mathbb{S}_{3,+}^{n-1}$ be an eigenvector to $\rho(H_i)$. By Lemma \ref{lm_Liu}, we denote $x_l=x_v$ for $v\in V_l^i$, where $l=1,2,3$. By Proposition \ref{Hou}, $\rho(H_i)=O(n^2)$ for $2\leq i\leq k+1$.
 Applying the eigenvector equation, we have
\begin{equation}\label{nontrivial eigenvector equation3}
\begin{cases}
\rho(H_i)x_1^{2}&=\binom{2i-3}{2}x_1^2+(2i-3)(3k-3i+5)x_1x_2+(2i-3)(n-3k-3+i)x_1x_3\\
&+\binom{3k-3i+5}{2}x_2^2+(3k-3i+5)(n-3k-3+i)x_2x_3+\binom{n-3k-3+i}{2}x_3^2\\
\rho(H_i)x_2^{2}&=\binom{2i-2}{2}x_1^2+(2i-2)(3k-3i+4)x_1x_2+(2i-2)(n-3k-3+i)x_1x_3\\
&+\binom{3k-3i+4}{2}x_2^2\\
\rho(H_i)x_3^{2}&=\binom{2i-2}{2}x_1^2+(2i-2)(3k-3i+5)x_1x_2+(2i-2)(n-3k-4+i)x_1x_3
\end{cases}
\end{equation}
Let $$E_1=\big\{e:e\in E(H_i), |e\cap [3k+3,n]|\geq 1\big\},$$ 
$$E_2=\left\{e\in \binom{[n]}{3}: |e\cap V_2^i|=2, |e\cap \left(V_3^i\setminus\{3k+3, \ldots, n\}\right)|=1\right\},$$
$$E_3=\left\{e\in \binom{[n]}{3}: |e\cap V_2^i|=1, |e\cap \left(V_3^i\setminus\{3k+3, \ldots,n\}\right)|=2\right\},$$
and
$$E_4=\left\{e\in \binom{[n]}{3}: |e\cap \left(V_3^i\setminus\{3k+3, \ldots, n\}\right)|=3\right\}.$$
For any $2\leq i\leq k+1$, it is easy to see that
$$E(H_1)=\big(E(H_i)\cup (E_2\cup E_3\cup E_4)\big)\backslash E_1.$$
Then, it follows from that
\begin{equation}\label{diff}
\begin{aligned}
&\rho(H_1)-\rho(H_i)\\
&\geq 3\left(\sum_{e\in E_2\cup E_3\cup E_4}\mathbf{x}(e)-\sum_{e\in E_1}\mathbf{x}(e)\right)\\
&=3\left[\binom{3k-3i+5}{2}(i-1)x_2^2x_3+(3k-3i+5)\binom{i-1}{2}x_2x_3^2+\binom{i-1}{3}x_3^3\right]\\
&-3(2i-2)(3k-3i+5)(n-3k-2)x_1x_2x_3-3\binom{2i-2}{2}(n-3k-2)x_1^2x_3\\
&-3(2i-2)(i-1)(n-3k-2)x_1x_3^2-3(2i-2)\binom{n-3k-2}{2}x_1x_3^2.
\end{aligned}
\end{equation}
 
Next, we distinguish the following two cases.
\begin{case}
  $i=\theta(n)$.
\end{case}

For $l=1,2,3$, let $x_l=c_ln^{\alpha_l}+o(n^{\alpha_l})$, where $c_l$ is a constant. By (\ref{nontrivial eigenvector equation3}), we can obtain $\alpha_1=\alpha_2=\alpha_3$. Hence, combining (\ref{diff}),
\begin{equation*}
\begin{aligned}
\rho(H_1)-\rho(H_i)&=\theta(n^3)n^{3\alpha_1}-\theta(n^2)n^{3\alpha_1}\\
&>0.
\end{aligned}
\end{equation*}

\begin{case}
  $i=o(n)$.
\end{case}
Let $i=c_in^{\alpha_i}+o(n^{\alpha_i})$ and $n-3k=\theta(n^{\beta})$, where $\alpha_i<1$ and $\beta<1$.
By (\ref{diff}), we derive
\begin{equation*}
\begin{aligned}
\rho(H_1)-\rho(H_i)&=\theta(n^{2+\alpha_i})n^{3\alpha_1}-\theta(n^{1+\alpha_i+\beta})n^{3\alpha_1}\\
&>0.
\end{aligned}
\end{equation*}
\c

Recall that $\rho(H')\leq \max\{\rho(H_1), \rho(H_2), \ldots, \rho(H_{k+1})\}$ and $\rho(H')\geq \rho(H_{1})$. Thus, by Claim \ref{spectral comparison2}, $$\rho(H')=\rho(H_{1})=\rho(\mathcal{A}^3_3).$$ This completes the proof.
\q

\subsection{Proof of the Theorem \ref{HMTT}}

%\begin{theorem}
%Let $n$ be a positive integer with $n\geq 8$. If $H$ is a non-trivial intersecting $3$-graph on $n$ vertex. Then
%$$\rho(H)\leq \max\{\rho(H_1),\rho(H_2)\},$$
%where $H_1=\left\{e\in \binom{[n]}{3}:  |e\cap \{1,2,3\}|\geq 2\right\}$  and $H_2=\left\{e\in \binom{[n]}{3}:  1\in e, |e\cap \{2,3,4\}|\geq 1\right\}\cup \{2,3,4\}$.
%\end{theorem}
\p
Let $H$ be a non-trivial intersecting $3$-graph on vertex set $[n]$ with the maximum spectral radius.  It is easy to see that $\mathcal{H}_1$ and $\mathcal{F}^3_1$ are two non-trivial intersecting $3$-graphs. So, 
\begin{equation}\label{lower bound}
\rho(H)\geq \max\{\rho(\mathcal{H}_1), \rho(\mathcal{F}^3_1)\}.
\end{equation}

By Lemma \ref{Matching}, we can  apply a number of shifting operations to $H$ while preserving the intersecting property of $H$. The only trouble that might occur in the shifting process is that $H$ is non-trivial but $H_{j\rightarrow i}$ is trivial. Hence, we firstly show the following claim.
\begin{claim}
  There is a shifted non-trivial intersecting $3$-graph with the maximum spectral radius.
\end{claim}
\p
Suppose that $H'$ is a non-trivial intersecting $3$-graph obtained by applying a number of shifting operations to $H$ but $H'_{j\rightarrow i}$ is trivial, where $i,j\in [n]$.  By renaming the vertices we may assume that $i=1$ and $j=2$ in $H'$. In this case, it follows that $e\cap \{1,2\}\neq \emptyset$ for any $e\in E(H')$, and there exist edges $e_1$ and $e_2$ such that $e_1\cap \{1,2\}=\{1\}$ and $e_2\cap \{1,2\}=\{2\}$. 

Let $H^{(1)}$ be a $3$-graph with vertex set $[n]$ and edge set $E(H^{(1)})=\{e\in \binom{[n]}{3}:\{1,2\}\subset e\}.$ 
We  show that $e\in E(H')$ for any $e\in E(H^{(1)})$. Otherwise, suppose that there exists an edge $e'\in E(H^{(1)})$ and $e'\notin E(H')$. Let $H''$ be a $3$-graph with $V(H'')=V(H')$ and $E(H'')=E(H')\cup\{e'\}$. Clearly, $\rho(H'')> \rho(H')$. Meanwhile, based on the fact that $H'$ is a non-trivial intersecting $3$-graph  and $e\cap \{1,2\}\neq \emptyset$ for any $e\in E(H')$, we have that $H''$ is also a non-trivial intersecting $3$-graph.  
By Lemma \ref{spectral} and the choice of $H$, we have $\rho(H')=\rho(H)$. Then, $\rho(H'')>\rho(H')=\rho(H)$, which contradicts the choice of $H$.

Note that $H^{(1)}_{j\rightarrow i}=H^{(1)}$ for all $3\leq i<j\leq n$. 
Then, starting from $H'$, if we apply all possible shifting operations $S_{j\rightarrow i}$ for $3\leq i<j\leq n$, then the resulting $3$-graph is non-trivial since it still contains $H^{(1)}$ and there exist edges $e_3$ and $e_4$ such that $e_3\cap \{1,2\}=\{1\}$ and $e_4\cap \{1,2\}=\{2\}$. So we may assume that we can obtain that $H'$ is shifted on $[3,n]$ , that is, if $e\in E(H')$ and $e\cap\{i,j\}=j$ then $(e\backslash \{j\})\cup \{i\}\in E(H')$ for all $3\leq i<j\leq n$. Hence, $\{2,3,4\}, \{1, 3, 4\}\in E(H')$. Recall that $\{1, 2, 3\}, \{1, 2, 4\}\in E(H')$.
Let $H^{(2)}$ be a $3$-graph with vertex set $[n]$ and edge set $E(H^{(2)})=\binom{[4]}{3}.$  Then, $H^{(2)}$ is a subgraph of $H'$ and $H^{(2)}$ is unchanged under any shifting $S_{j\rightarrow i}$, $1\leq i<j\leq n$. From now on we
never create any trivial intersecting families in the shifting process, and  in the end we get a shifted non-trivial intersecting $3$-graph.
\c

Let $H'$ be a shifted non-trivial intersecting $3$-graph obtained by applying a number of shifting operations to $H$. Then, by Proposition \ref{shiftedP}, $\{2,3,4\}\in E(H')$. Meanwhile, by the choice of $H$ and Lemma \ref{spectral}, $\rho(H)=\rho(H')$. 
Let $e_1=\{1,4, 6\}$ and $e_2=\{2,3,5\}$.  Then $e_1$ and $e_2$ cannot both be in $E(H')$. Next, we distinguish the following two cases.

\begin{case}
  $e_1=\{1,4,6\}\notin E(H')$.
\end{case}

Based on the fact that $H'$ is a non-trivial intersecting $3$-graph with the maximum spectral radius, we show that $\{1,4,5\}\notin E(H')$. We prove this by contradiction. Suppose $\{1,4,5\}\in E(H')$. Note that there are the following two subcases:

\begin{case2}
$\{2,4,5\}\in E(H')$.
\end{case2}
%\textbf{Case 1.1.} $\{2,4,5\}\in E(H')$.

Since $\mathcal{H}_1$ is a non-trivial intersecting $3$-graph, we have that  $\rho(H')\geq \rho(\mathcal{H}_1)$. 
From the eigenvector equation (\ref{r eigenvector equation}), we obtain $\rho(\mathcal{H}_1)(\rho(\mathcal{H}_1)-1)^2=12(n-3)^2.$  Then, for sufficiently large $n$,  $\rho(H')\geq \rho(\mathcal{H}_1)>6$.
We now show that $\{3,4,5\}\notin E(H')$. Suppose for contradiction that $\{3,4,5\}\in E(H')$. Since $\{1,2,c\}\cap \{3,4,5\}=\emptyset$ for every $c\geq 6$, it follows that $\{1,2,c\}\notin E(H')$.
Recall that $H'$ is shifted.  By Proposition \ref{shiftedP}, we have $d_{H'}(c)=0$ for any $c\geq 6$ and $E(H')= \binom{[5]}{3}$. Hence, $\rho(H')=6$,  a contradiction.  
By Proposition \ref{shiftedP}, it is easy to see that $\{a,b,c\}\in E(H')$ satisfies $a\leq 2$. Note that $\{1,4,6\}\notin E(H')$.
Then, for any $\{a,b,c\}\in E(H')$, $b\leq 3$ or $c\leq 5$. Consequently, $$E(H')\subseteq\{e\in \binom{[n]}{3}: |e\cap\{1,2,3\}|\geq 2\}\cup \big\{\{2,4,5\},\{1,4,5\}\big\}.$$  Furthermore, since $\{2,4,5\}\in E(H')$ and $H'$ is a intersecting $3$-graph, we get that $\{1,3,c\}\notin E(H')$ for $c\geq 6$. Let
$$E'=\big\{\{1,3,4\},\{2,3,4\},\{1,3,5\},\{2,3,5\},\{1,4,5\},\{2,4,5\}\big\}.$$
Then, $$E(H')=\left\{e\in \binom{[n]}{3}: \{1,2\}\subseteq e\right\}\cup E'.$$

Next, to complete the proof, we derive a contradiction to the choice of $H'$ by showing that  $\rho(H')<\rho(\mathcal{F}^3_1)$. Let $V_1(H')=\{1,2\}$, $V_2(H')=\{3,4,5\}$  and $V_3(H')=\{6,\ldots, n\}$. For any $i,j\in V_l(H')$, it is easy to see that $i$ and $j$ are equivalent in $H'$, where $l=1,2,3$.
Let $\mathbf{x}\in \mathbb{S}_{3,+}^{n-1}$ be an eigenvector to $\rho(H')$. By Lemma \ref{lm_Liu}, we denote $x_l=x_v$ for $v\in V_l(H')$, where $l=1,2,3$.  Applying the eigenvector equation, we have
\begin{equation}\label{nontrivial eigenvector equation}
\begin{cases}
\rho(H')x_1^{2}&=(n-5)x_1x_3+3x_1x_2+3x_2^2\\
\rho(H')x_2^{2}&=x_1^2+4x_1x_2\\
\rho(H')x_3^{2}&=x_1^2
\end{cases}
\end{equation}
Combining $(\ref{nontrivial eigenvector equation})$,  we have
\begin{equation*}
\begin{aligned}
\rho(H')(x_1^{2}-x_2^{2})&=(n-5)x_1x_3+3x_1x_2+3x_2^2-x_1^2-4x_1x_2\\
&=(n-5)x_1x_3+x^2_2+x_2(x_2-x_1)+x^2_2-x_1^2,
\end{aligned}
\end{equation*}
which implies that
\begin{equation*}
\begin{aligned}
-\rho(H')(x_1+x_2)=2x_2+x_1+\frac{(n-5)x_1x_3+x_2^2}{(x_2-x_1)}.
\end{aligned}
\end{equation*}
Since $\rho(H')(x_1+x_2)>0, 2x_2+x_1>0$ and $(n-5)x_1x_3+x_2^2>0$, it follows that $x_2<x_1$. Furthermore, we claim that $(n-5)x_3>x_2$. Suppose, for contradiction, that $(n-5)x_3\leq x_2$. Combining $(\ref{nontrivial eigenvector equation})$ and $x_2<x_1$, we get
\begin{equation*}
\begin{aligned}
(n-5)^2x_1^2=(n-5)^2\rho x^2_3\leq \rho x^2_2=x_1^2+4x_1x_2<5x_1^2,
\end{aligned}
\end{equation*}
a contradiction.
Let $E_1=E(H')\backslash E(\mathcal{F}^3_1)=\{2,4,5\}$ and $$E_2=E(\mathcal{F}^3_1)\backslash E(H')=\left\{\{1,3,6\}, \ldots, \{1,3,n\}\right\}\cup\left\{\{1,4,6\}, \ldots, \{1,4,n\}\right\}.$$ Then, 
\begin{equation}\label{nontrivial spectral diff}
\begin{aligned}
\rho(\mathcal{F}^3_1)-\rho(H')&\geq 3\left(\sum_{e\in E_2}\mathbf{x}(e)-\sum_{e\in E_1}\mathbf{x}(e)\right)\\
&=3\left[2(n-5)x_1x_2x_3-x_1x_2^2\right]\\
&=3x_1x_2\big(2(n-5)x_3-x_2\big)\\
&>0,
\end{aligned}
\end{equation}
where the last inequality follows from that fact that $(n-5)x_3>x_2$.

\begin{case2}
$\{2,4,5\}\notin E(H')$.
\end{case2}
%\textbf{Case 1.2.} $\{2,4,5\}\notin E(H')$.

Since $\{1,4,5\}\in E(H')$, $\{1,4,6\}\notin E(H')$, $\{2,4,5\}\notin E(H')$ and $H'$ is a non-trivial intersecting shifted $3$-graph with the maximum spectral radius, by Proposition \ref{shiftedP}, we have 
$$E(H')=\left\{e\in \binom{[n]}{3}: \{1,2\}\subseteq e\ \text{or}\ \{1,3\}\subseteq e\right\}\cup \big\{\{2,3,4\},\{2,3,5\},\{1,4,5\}\big\}.$$

Next, to complete the proof, we derive a contradiction to the choice of $H'$ by showing that  $\rho(H')<\rho(\mathcal{H}_1)$ or $\rho(H')<\rho(\mathcal{F}^3_1)$.
Let $V_1(H')=\{1\}$, $V_2(H')=\{2,3\}$, $V_3(H')=\{4,5\}$ and $V_4(H')=\{6,\ldots, n\}$. For any $i,j\in V_l(H')$, it is easy to see that $i$ and $j$ are equivalent in $H'$, where $l=2,3,4$.
Let $\mathbf{x}$ be an eigenvector to $\rho(H')$. Using Lemma \ref{lm_Liu}, we denote $x_l=x_v$ for $v\in V_l(H')$, where $l=2,3,4$.  Form eigenvector equation, we have
\begin{equation}\label{nontrivial eigenvector equation2}
\begin{cases}
\rho(H')x_1^{2}&=2(n-5)x_2x_4+4x_2x_3+x_2^2+x^2_3\\
\rho(H')x_2^{2}&=(n-5)x_1x_4+x_1x_2+2x_1x_3+2x_2x_3\\
\rho(H')x_3^{2}&=x_2^2+x_1x_3+2x_1x_2\\
\rho(H')x_4^{2}&=2x_1x_2
\end{cases}
\end{equation}
Using $(\ref{nontrivial eigenvector equation2})$, we  obtain
\begin{equation*}
\begin{cases}
\begin{aligned}
\rho(H')(x_1^{2}-x_2^2)&=(n-5)x_2x_4+x_3^2+(x_2-x_1)\left((n-5)x_4+x_2+2x_3\right)\\
\rho(H')(x_2^{2}-x_3^2)&=(n-5)x_1x_4+x_1(x_3-x_2)+x_2(x_3-x_2)+x_2x_3,
\end{aligned}
\end{cases}
\end{equation*}
which can be rearranged as
\begin{equation*}
\begin{cases}
-\rho(H')(x_1+x_2)&=\frac{(n-5)x_2x_4+x_3^2}{x_2-x_1}+(n-5)x_4+x_2+2x_3\\
-\rho(H')(x_2+x_3)&=\frac{(n-5)x_1x_4+x_2x_3}{x_3-x_2}+x_1+x_2
\end{cases}
\end{equation*}
Since $\rho(H')$ and $x_i>0$ for $i=1,2,3,4$, we conclude that $x_1>x_2>x_3$. Furthermore, we show that $(n-5)x_4>x_3$. Indeed, suppose for contradiction that $(n-5)x_4\leq x_3$. Combining (\ref{nontrivial eigenvector equation2}), we have 
\begin{equation*}
\begin{aligned}
2(n-5)^2x_1x_2=(n-5)^2\rho x^2_4\leq \rho x^2_3=2x_1x_2+x_1x_3+x_2^2<4x_1x_2,
\end{aligned}
\end{equation*}
 a contradiction. 
 
We now analyze the two cases according to the relation between $x_2^2$ and $x_1x_3$.
First, if $x_2^2>x_1x_3$, then
\begin{equation*}
\begin{aligned}
\rho(\mathcal{H}_1)-\rho(H')&\geq 3\left(\sum_{e\in E(\mathcal{H}_1)\backslash E(H')}\mathbf{x}(e)-\sum_{e\in E(H')\backslash E(\mathcal{H}_1)}\mathbf{x}(e)\right)\\
&=3\left[(n-5)x_4x_2^2-x_1x_3^2\right]\\
&>0,
\end{aligned}
\end{equation*}
contradicting the choice of $H'$.  On the other hand, if $x_2^2\leq x_1x_3$, then
\begin{equation*}
\begin{aligned}
\rho(\mathcal{F}^3_1)-\rho(H')&\geq 3\left(\sum_{e\in E(\mathcal{F}^3_1)\backslash E(H')}\mathbf{x}(e)-\sum_{e\in E(H')\backslash E(\mathcal{F}^3_1)}\mathbf{x}(e)\right)\\
&=3\left[(n-5)x_4x_1x_3-x_3x_2^2\right]\\
&>0,
\end{aligned}
\end{equation*}
which contradicts the choice of $H'$.  

Accordingly, we know that $\{1,4,5\}, \{1,4,6\}\notin E(H')$ and $\{2,3,4\}\in E(H')$.  Applying Proposition \ref{shiftedP}, we conclude that 
$$E(H')=\left\{e\in \binom{[n]}{3}:  |e\cap \{1,2,3\}|\geq 2\right\}=E(\mathcal{H}_1).$$

\begin{case}
  $e_2=\{2,3,5\}\notin E(H')$.
\end{case}

Since $\{2,3,4\}\in E(H')$ and $\{2,3,4\}\cap \{1,5,6\}=\emptyset$, it follows that $\{1,5,6\}\notin E(H')$.
Based on the face that $H'$ is a non-trivial intersecting shifted $3$-graph with the maximum spectral radius and $\{1,5,6\}, \{2,3,5\} \notin E(H')$, by Proposition \ref{shiftedP}, we deduce that
$$E(H')=\left\{e\in \binom{[n]}{3}:  1\in e, |e\cap \{2,3,4\}|\geq 1\right\}\cup \big\{\{2,3,4\}\big\}=E(\mathcal{F}^3_1).$$

Through the analysis of the above two cases and (\ref{lower bound}), we obtain $$\rho(H)=\rho(H')= \max\{\rho(\mathcal{H}_1), \rho(\mathcal{F}^3_1)\}.$$
In order to complete the proof, we show that $\rho(\mathcal{H}_1)>\rho(\mathcal{F}^3_1)$ for sufficiently large  $n$.

Let $V_1(\mathcal{H}_1)=\{1,2,3\}$, $V_2(\mathcal{H}_1)=\{4,\ldots, n\}$. For any $i,j\in V_l(\mathcal{H}_1)$, it is easy to see that $i$ and $j$ are equivalent in $\mathcal{H}_1$, where $l=1,2$.
Let $\mathbf{x}\in \mathbb{S}_{r,+}^{n-1}$ be an eigenvector to $\rho(\mathcal{H}_1)$. By Lemma \ref{lm_Liu}, we denote $x_l=x_v$ for $v\in V_l(\mathcal{H}_1)$, where $l=1,2$.  Applying the eigenvector equation, we have
\begin{equation*}
\begin{cases}
\rho(\mathcal{H}_1)x_1^{2}&=2(n-3)x_1x_2+x_1^2\\
\rho(\mathcal{H}_1)x_2^{2}&=3x_1^2,
\end{cases}
\end{equation*}
which implies that
\begin{equation}\label{eqSpectral1}
(\rho(\mathcal{H}_1)-1)^2\rho(\mathcal{H}_1)=12(n-3)^2.
\end{equation}
Assume $\rho(\mathcal{H}_1)=cn^{\alpha}+o(n^{\alpha})$, where $c$ is a positive constant. Then, combining (\ref{eqSpectral1}), it follows that
\begin{equation*}
c^3n^{3\alpha}+o(n^{3\alpha})=12(n-3)^2.
\end{equation*}
Comparing the leading orders of both sides, we deduce that $c=\sqrt[3]{12}$ and $\alpha=\frac{2}{3}$. Thus, $\rho(\mathcal{H}_1)=\sqrt[3]{12}n^{\frac{2}{3}}+o(n^{\frac{2}{3}})$.

We partition $V(\mathcal{F}^3_1)$ into three disjoint subsets $V_1(\mathcal{F}^3_1)=\{1\}$, $V_2(\mathcal{F}^3_1)=\{2, 3, 4\}$ and $V_3(\mathcal{F}^3_1)=\{5,\ldots, n\}$. For any $i,j\in V_l(\mathcal{F}^3_1)$, it is easy to see that $i$ and $j$ are equivalent in $\mathcal{F}^3_1$, where $l=1,2,3$.
Let $\mathbf{x}\in \mathbb{S}_{3,+}^{n-1}$ be an eigenvector to $\rho(\mathcal{F}^3_1)$. By Lemma \ref{lm_Liu}, we denote $x_l=x_v$ for $v\in V_l(\mathcal{F}^3_1)$, where $l=1,2,3$.  Applying the eigenvector equation, we obtain
\begin{equation}\label{H2}
\begin{cases}
\rho(\mathcal{F}^3_1)x_1^{2}&=3(n-4)x_2x_3+3x_2^2\\
\rho(\mathcal{F}^3_1)x_2^{2}&=x_2^2+2x_2x_1+(n-4)x_1x_3\\
\rho(\mathcal{F}^3_1)x_3^{2}&=3x_1x_2.
\end{cases}
\end{equation}
Setting $u=\frac{x_1}{x_2}$ and $t=\frac{x_3}{x_2}$, we simplify  (\ref{H2}) to
\begin{equation}\label{SEH2}
\begin{cases}
\rho(\mathcal{F}^3_1)u^2&=3(n-4)t+3\\
\rho(\mathcal{F}^3_1)&=1+2u+(n-4)ut\\
\rho(\mathcal{F}^3_1)t^2&=3u.
\end{cases}
\end{equation}
Solving for $t$ from the third equation in (\ref{SEH2}), we obtain
\begin{equation*}
t=\sqrt{\frac{3u}{\rho(\mathcal{F}^3_1)}}.
\end{equation*}
Substituting this into the first two equations in (\ref{SEH2}), we have 
\begin{equation}\label{SH2}
\begin{cases}
\rho(\mathcal{F}^3_1)u^2&=3(n-4)\sqrt{\frac{3u}{\rho(\mathcal{F}^3_1)}}+3\\
\rho(\mathcal{F}^3_1)&=1+2u+(n-4)u\sqrt{\frac{3u}{\rho(\mathcal{F}^3_1)}}.
\end{cases}
\end{equation}
Let $\rho(\mathcal{F}^3_1)=cn^{\beta}+o(n^{\beta})$ and $u=c_1n^{\beta_1}+o(n^{\beta_1})$, where $c$ and $c_1$ are positive constants. Then, substituting these asymptotic forms into (\ref{SH2}) and retaining only the leading terms, we derive
\begin{equation*}
\begin{cases}
c^{\frac{3}{2}}c_1^2n^{\frac{3\beta}{2}+2\beta_1}+o(n^{\frac{3\beta}{2}+2\beta_1})&=
3n\sqrt{3c_1n^{\beta_1}}+o(n^{1+\frac{\beta_1}{2}})\\
c^{\frac{3}{2}}n^{\frac{3\beta}{2}}+o(n^{\frac{3\beta}{2}})&=c_1n^{1+\beta_1}\sqrt{3c_1n^{\beta_1}}
+o(n^{1+\frac{3\beta_1}{2}}),
\end{cases}
\end{equation*}
which simplifies to
\begin{equation}\label{SSH2}
\begin{cases}
c^{\frac{3}{2}}n^{\frac{3\beta}{2}}+o(n^{\frac{3\beta}{2}})&=
\frac{3n\sqrt{3c_1n^{\beta_1}}}{c_1^2n^{2\beta_1}}
+o(n^{1-\frac{3\beta_1}{2}})\\
c^{\frac{3}{2}}n^{\frac{3\beta}{2}}+o(n^{\frac{3\beta}{2}})&=c_1n^{1+\beta_1}\sqrt{3c_1n^{\beta_1}}
+o(n^{1+\frac{3\beta_1}{2}}).
\end{cases}
\end{equation}
By equating the leading terms of the two equations in (\ref{SSH2}), we obtain
$$\frac{3n\sqrt{3c_1n^{\beta_1}}}{c_1^2n^{2\beta_1}}=c_1n^{1+\beta_1}\sqrt{3c_1n^{\beta_1}},$$
which gives $c_1=\sqrt[3]{3}$ and $\beta_1=0$.
Substituting these values into the second equations in (\ref{SSH2}), we obtain that $c=3^{\frac{2}{3}}$ and $\beta=\frac{2}{3}$. Hence, $\rho(\mathcal{F}^3_1)=\sqrt[3]{9}n^{\frac{2}{3}}+o(n^{\frac{2}{3}})$.

Therefore, we conclude that $\rho(\mathcal{H}_1)=\sqrt[3]{12}n^{\frac{2}{3}}+o(n^{\frac{2}{3}})$ and $\rho(\mathcal{F}^3_1)=\sqrt[3]{9}n^{\frac{2}{3}}+o(n^{\frac{2}{3}})$. Since $\sqrt[3]{12}>\sqrt[3]{9}$, it follows that  $\rho(\mathcal{H}_1)>\rho(\mathcal{F}^3_1)$ for sufficiently large $n$.
\q

\end{document}